\documentclass[10pt]{amsart}
\usepackage{mathptmx}
\usepackage{amsmath}     
\usepackage{amssymb}
\usepackage{array}
\usepackage{geometry}
\usepackage[bookmarks=true,colorlinks=true, pdfstartview=FitV, linkcolor=black, citecolor=blue, urlcolor=black]{hyperref}
\usepackage{movie15}

\usepackage{color}
\definecolor{DarkRed}{rgb}{0.55,.00,0.2}
\definecolor{DarkGrey}{rgb}{0.35,.35,0.35}

\theoremstyle{definition}

\theoremstyle{remark}

\numberwithin{equation}{section}

\begin{document}

\title{A method of  composition orthogonality and  new sequences of orthogonal  polynomials and functions for non-classical weights}

\author{{\bf S. Yakubovich} \\
 {\em {Department of Mathematics, Fac. Sciences of University of Porto,\\ Rua do Campo Alegre,  687; 4169-007 Porto (Portugal) }}}
  \thanks{ E-mail: syakubov@fc.up.pt} 
\vspace{3mm}

\thanks{ The work was  partially supported by CMUP, which is financed by national funds through FCT(Portugal),  under the project with reference UIDB/00144/2020. }


\date{\today}



\begin{abstract}  A new method of composition orthogonality is introduced. It is applied to generate new sequences of orthogonal polynomials and functions. In particular, classical orthogonal polynomials are interpreted in the sense of composition orthogonality. Finally, new sequences of  orthogonal polynomials  with respect to the weight function $\rho^2_\nu(x), \rho_{\nu}(x)= 2 x^{\nu/2} K_\nu(2\sqrt x),\ x >0,  \nu \ge 0$, where  $K_\nu(z)$ is the modified Bessel function or Macdonald function, are investigated.  Differential properties,  recurrence relations, explicit representations, generating functions and Rodrigues-type formulae  are obtained.  The corresponding multiple  orthogonal polynomials  are exhibited. \\
\end{abstract}

\maketitle

\markboth{\rm \centerline{ S. Yakubovich}}{Composition orthogonality and new sequences of orthogonal  polynomials and functions}

\section{Introduction and preliminary results}

In studying Prudnikov's  sequence of orthogonal polynomials  (see \cite{YAP})  with  the weight function $x^\alpha \rho_\nu(x)$, where \  $ \rho_{\nu}(x)= 2 x^{\nu/2} K_\nu(2\sqrt x),\ x >0, \nu \ge 0,\ \alpha > -1$ and  $K_\nu(z)$ is  the Macdonald or modified Bessel function  \cite{Bateman}, Vol. II the author interpreted this sequence in terms of the  Laguerre composition orthogonality, involving the differential operator $\theta= xDx$, where $D= d/dx$.  The main aim of this paper is to extend the method, employing classical orthogonal polynomials (Hermite, Laguerre, Jacobi) \cite{Bateman} to find new sequences of orthogonal polynomials with non-classical weights such as, for instance,  the square of Macdonald function and other hypergeometric functions.  In fact,   we will define the composition orthogonality in the following way.

{\bf Definition 1.} {\it Let  $\omega(t), \varphi(t) \ t \in [a,b],\  -\infty \le a < b \le \infty$ be  nonnegative  functions, $\theta = tDt$.    Let $f, g$ be complex-valued functions such that the product of two operators $f(\theta) g(\theta)$ commutes. Then  $f,\ g $ are compositionally orthogonal with respect to the measure $\omega(t)dt$ relatively to the function $\varphi$ if }
$$\int_a^b  f(\theta) g(\theta) \{ \varphi(t)\} \omega(t) dt = 0.\eqno(1.1)$$ 
We will see in the sequel that this definition suits well with the vector space of polynomials over $\mathbb{R}$, when  the composition orthonormality of the sequence $\{P_n\}_{n\ge 0}$ is defined by

$$\int_a^b  P_n(\theta) P_m(\theta) \{ \varphi(t)\} \omega(t) dt = \delta_{n,m},\eqno(1.2)$$ 
where $\delta_{n,m},\   n,m\in\mathbb{N}_{0}$ is the Kronecker symbol.   Moreover, for some class of functions $\varphi$ it is possible to transform the left-hand side of the equality (1.2)  to the usual  orthogonality  with respect to a new weight function.   Then  since the $n$-th power of the operator $\theta$ satisfies the Viskov-type identity \cite{Viskov}

$$ \theta^n = \left( xDx\right)^n = x^n D^n x^n,\quad   n \in\mathbb{N}_{0},\eqno(1.3)$$
we can get via integration by parts new properties of the sequence $\{P_n\}_{n\ge 0}$ and its relationship, for instance, with classical orthogonal polynomials.  In fact, let us consider a class of functions $\varphi$ representable by the modified Laplace transform of some nonnegative function $\psi$

$$\varphi(t)= {1\over t} \int_0^\infty e^{-x/t} \psi(x) dx,\ t \in [a,b] \subset (0, \infty).\eqno(1.4)$$
Hence it is easily seen that 

$$\theta^k \left\{ t^{-1} e^{-x/t} \right\}  = \left( t D t\right)^k \left\{ t^{-1} e^{-x/t} \right\} = x^k  t^{-1}  e^{-x/t},\quad k \in \mathbb{N}_0,\eqno(1.5)$$
and therefore, differentiating under the integral sign in (1.4), we derive

$$ \theta^k \left\{ \varphi(t)  \right\} = {1\over t} \int_0^\infty e^{-x/t}  x^k \psi(x) dx,\quad  k \in \mathbb{N}_0. \eqno(1.6)$$
It is indeed allowed, for instance,  due to the assumed convergence of the integral 

$$ \int_0^\infty e^{-x/b}  x^k \psi(x) dx < \infty, \quad   k \in \mathbb{N}_0.\eqno(1.7)$$
Consequently, returning to (1.2) we write its left-hand side in the form

$$ \int_a^b  P_n(\theta) P_m(\theta) \{ \varphi(t)\} \omega(t) dt  = \int_a^b   \int_0^\infty e^{-x/t} P_n(x) P_m(x) \psi(x)  \omega(t) {dx dt \over t}.\eqno(1.8) $$ 
The interchange of the order of integration on the right-hand side in (1.8) is permitted by Fubini's theorem under the imposed condition

$$ \int_a^b   \int_0^\infty e^{-x/t}  x^k \psi(x)  \omega(t) {dx dt \over t} < \infty,\quad  k \in \mathbb{N}_0, \eqno(1.9)$$
and, combining with (1.2), we find the equalities

$$\int_a^b  P_n(\theta) P_m(\theta) \{ \varphi(t)\} \omega(t) dt =  \int_0^\infty  P_n(x) P_m(x) \psi(x) \Omega (x) dx = \delta_{n,m},\eqno(1.10)$$ 
where 

$$ \Omega (x) = \int_a^b  e^{-x/t}  \omega(t) {dt \over t},\quad  x >0.\eqno(1.11) $$
Thus we see that the sequence  $\{P_n\}_{n\ge 0}$ is orthonormal over $(0,\infty)$ with respect to the measure $\Omega (x) dx$.  Moreover, one can consider 
the left-hand side of the first equality in (1.10) as an  inner product on the vector space of  polynomials over $\mathbb{R}$  (a pre-Hilbert space)

$$\langle p, q \rangle =  \int_a^b  p (\theta) q (\theta) \{ \varphi(t)\} \omega(t) dt,\eqno(1.12)$$
inducing the norm by the equality

$$ ||p||= \sqrt{\langle p, p \rangle} = \left( \int_0^\infty  p^2(x) \psi(x) \Omega (x) dx\right)^{1/2}.\eqno(1.13)$$

\section{The use of classical orthogonal polynomials}

\subsection{Laguerre polynomials}
We begin to consider sequences of orthogonal polynomials which are compositionally orthogonal   with respect to the measure $t^\nu e^{-t} dt$ over $\mathbb{R}_+$ related to   Laguerre polynomials $\{L_n^\nu\}_{n\ge 0},\ \nu > -1$. In fact, letting $\omega(t)=  t^\nu e^{-t},\ t > 0$ in (1.2), we get

$$\int_0^\infty  P_n(\theta) P_m(\theta) \{ \varphi(t)\}  t^\nu e^{-t} dt = \delta_{n,m}.\eqno(2.1)$$ 
The corresponding integral (1.11) is calculated in \cite{Bateman}, Vol. II, and we obtain

$$\Omega(x)= \int_0^\infty  e^{-x/t -t}  t^{\nu-1} dt = 2 x^{\nu/2} K_\nu\left( 2\sqrt x\right) \equiv \rho_\nu(x),\ x > 0,\eqno(2.2)$$
where $K_\nu(z)$ is the modified Bessel function or Macdonald function \cite{YaL}.  The  function  $\rho_{\nu}$  has  the Mellin-Barnes integral representation in the form (cf. \cite{YAP})
$$
\rho_\nu(x)=  \frac{1}{2\pi i} \int_{\gamma-i\infty}^{\gamma+i\infty} \Gamma(\nu+s) \Gamma (s) x^{-s} ds\ , \quad x, \gamma \in \mathbb{R}_{+},\ \nu  \in \mathbb{R},\eqno(2.3)
$$
where $\Gamma(z)$ is Euler's gamma-function \cite{Bateman}, Vol. I. The asymptotic behavior of the modified Bessel function at infinity and near the origin \cite{Bateman}, Vol. II gives the corresponding values for the  function $\rho_\nu,\ \nu \in \mathbb{R}$.  Precisely, we have
$$\rho_\nu (x)= O\left( x^{(\nu-|\nu|)/2}\right),\  x \to 0,\ \nu\neq 0, \quad  \rho_0(x)= O( \log x),\ x \to 0,\eqno(2.4)$$

$$ \rho_\nu(x)= O\left( x^{\nu/2- 1/4} e^{- 2\sqrt x} \right),\ x \to +\infty.\eqno(2.5)$$
Therefore, if the condition (cf. (1.9))

$$ \int_0^\infty   x^k \rho_\nu(x)  \psi(x) dx   < \infty,\quad  k \in \mathbb{N}_0 \eqno(2.6)$$
holds valid,  we arrive at the following proposition. 

{\bf Proposition 1.} {\it Let $\nu > -1, \varphi, \psi$  be  nonnegative functions defined on  $\mathbb{R}_+$ which are related by the modified Laplace transform $(1.4)$.  Then under condition $(2.6)$  the finiteness of the integral 

$$ \int_0^\infty   x^k e^{-x/M}  \psi(x) dx   < \infty,\quad  k \in \mathbb{N}_0 \eqno(2.7)$$
for some $M >0$ the sequence $\{P_n\}_{n\ge 0}$ of orthogonal polynomials with respect to the measure $\rho_\nu(x)  \psi(x) dx$ over $\mathbb{R}_+$ is compositionally orthogonal in the sense of Laguerre relatively to the function $\varphi$, i.e.}

$$  \int_0^\infty  P_n(\theta) P_m(\theta) \{ \varphi(t)\}  t^\nu e^{-t} dt= \int_0^\infty  P_n(x) P_m(x)  \rho_\nu(x)  \psi(x) dx = \delta_{n,m}.\eqno(2.8)$$

\begin{proof}  Since  via (2.7) the integral 

$$  {1\over t} \int_0^\infty e^{-x/t}  x^k \psi(x) dx,  \quad  k \in \mathbb{N}_0$$
converges uniformly with respect to  $t \in [1/M, M], M >0$, the consecutive differentiation under the integral sign is allowed, and we derive, recalling (1.6)

$$P_n(\theta) P_m(\theta) \{ \varphi(t)\}  =  P_n(\theta) P_m(\theta) \left\{ {1\over t} \int_0^\infty e^{-x/t} \psi(x) dx \right\}$$

$$ =   {1\over t} \int_0^\infty e^{-x/t}  P_n(x) P_m(x) \psi(x) dx.$$
Then

$$ \int_0^\infty  P_n(\theta) P_m(\theta) \{ \varphi(t)\}  t^\nu e^{-t} dt = \lim_{M\to \infty}   \int_{1/M}^M \int_0^\infty e^{-x/t}  P_n(x) P_m(x) \psi(x)   t^{\nu-1}  e^{-t} dx dt $$

$$=  \int_0^\infty  P_n(x) P_m(x)  \rho_\nu(x)  \psi(x) dx = \delta_{n,m},$$
where the latter equality is guaranteed by  condition (2.6) and Fubini's theorem. 

\end{proof}

{\bf Remark 1.} {\it Letting $\nu \ge 0, \psi(x)= x^\alpha,\ \alpha > -1$, we find the sequence $\{P^{\nu,\alpha}_n\}_{n\ge 0}$ of Prudnikov's orthogonal polynomials  studied in \cite{YAP}}

$$ \int_0^\infty  P^{\nu,\alpha}_n(x) P^{\nu,\alpha}_m(x)  \rho_\nu(x)  x^\alpha dx =  \Gamma(\alpha +1 ) \int_0^\infty  P_n(\theta) P_m(\theta) \{ t^\alpha\}  t^\nu e^{-t} dt = \delta_{n,m}.\eqno(2.9)$$

\subsection{Hermite polynomials}   Other interesting case is the Hermite orthogonality with respect to the measure $e^{-t^2} dt$ over $\mathbb{R}$. Our method will work on the following even extension  of the modified Laplace transform   (1.4)

$$\varphi( |t |)= {1\over| t|} \int_0^\infty e^{-x/|t| } \psi(x) dx,\ t \in \mathbb{R} \backslash{\{0\}}.\eqno(2.10)$$
Hence when  the Parseval identity for the Mellin transform \cite{Tit} suggests the equality

$$ \int_0^\infty  e^{-x/t -t^2}  {dt\over t} = {1\over 4\pi i} \int_{\gamma-i\infty}^{\gamma+i\infty} \Gamma\left({s\over 2}\right) \Gamma (s)\  x^{-s} ds,\quad   x,\gamma > 0. \eqno(2.11)$$
Invoking the duplication formula for the gamma function \cite{Bateman}, we get the right-hand side of (2.11) in the form

$$  {1\over 4\pi i} \int_{\gamma-i\infty}^{\gamma+i\infty} \Gamma\left({s\over 2}\right) \Gamma (s)\  x^{-s} ds = {1\over 4\pi^{3/2} i} \int_{\gamma/2 -i\infty}^{\gamma /2 +i\infty} \Gamma^2 \left(s\right) \Gamma \left({1\over 2}+ s \right)\ \left( {x\over 2}\right)^{- 2s} ds$$

$$  = {1\over 2\sqrt\pi}\  \rho_{1/2, 2} \left({x^2\over 4}\right),\eqno(2.12)$$
where by $\rho_{\nu, k}(x)$ we denote the ultra-exponential weight function introduced in \cite{YAP}

$$\rho_{\nu, k}(x) =  {1\over 2\pi i} \int_{\gamma-i\infty}^{\gamma+i\infty} \Gamma^k \left(s\right) \Gamma \left(\nu + s \right) x^{-s}ds, \  x,\nu,\gamma >0,  \ k \in \mathbb{N}_0 .\eqno(2.13)$$

{\bf Proposition 2.} {\it Let $t >0,\  x \in \mathbb{R},\  \varphi(t) , \psi (x)$  be  nonnegative functions  which are related by $(2.10)$.  If $\psi$ is even then under conditions $(2.7)$ and  

$$ \int_{-\infty}^\infty   |x|^k  \rho_{1/2, 2} \left({x^2\over 4}\right)  \psi(x) dx   < \infty,\quad  k \in \mathbb{N}_0 \eqno(2.14)$$
the sequence $\{P_n\}_{n\ge 0}$ of orthogonal polynomials with respect to the measure $2^{-1}\pi^{-1/2} \rho_{1/2, 2} \left({x^2/4}\right)  \psi(x) dx$ over $\mathbb{R}$ is compositionally orthogonal in the sense of Hermite  relatively to the function $\varphi(|t|)$, i.e.}

$$ \int_{-\infty}^\infty  P_n(\theta) P_m(\theta) \{ \varphi(|t|)\}   e^{-t^2} dt= \int_{-\infty}^\infty  P_n(x) P_m(x)   \rho_{1/2, 2} \left({x^2\over 4}\right)   \psi(x)  { dx\over \sqrt \pi}   = \delta_{n,m}.\eqno(2.15)$$

\begin{proof}  Indeed,  since

$$  \int_{-\infty}^\infty  P_n(\theta) P_m(\theta) \{ \varphi(|t|)\}   e^{-t^2} dt =  \int_{0}^\infty  \left[ P_n(\theta) P_m(\theta) + P_n(- \theta) P_m(- \theta)\right] \{ \varphi(t)\}   e^{-t^2} dt $$
we have due to (1.4), (1.6), (2.7)  

$$\left[ P_n(\theta) P_m(\theta) + P_n(- \theta) P_m(- \theta)\right] \{ \varphi(t)\} =  {1\over t} \int_0^\infty e^{-x/t}  \left[  P_n(x) P_m(x)  + P_n(- x) P_m(- x)  \right]  \psi(x) dx.\eqno(2.16)$$
Thus appealing to  (2.11), (2.12),  (2.14) and Fubini's theorem, we obtain finally  from (2.16)

$$   \int_{0}^\infty  \int_0^\infty   e^{-x/t- t^2} \left[  P_n(x) P_m(x)  + P_n(- x) P_m(- x)  \right]  \psi(x)  { dx dt \over t}$$

$$=  \int_0^\infty \left[  P_n(x) P_m(x)  + P_n(- x) P_m(- x)  \right]  \rho_{1/2, 2} \left({x^2\over 4}\right)   \psi(x) {dx\over \sqrt \pi}$$

$$= \int_{-\infty}^\infty  P_n(x) P_m(x)   \rho_{1/2, 2} \left({x^2\over 4}\right)   \psi(x)  { dx\over 2\sqrt \pi}   = \delta_{n,m}.$$

\end{proof}

\subsection{Jacobi  polynomials} The modified sequence $\{P_n^{\alpha,\beta}(2t-1)\}_{n\ge 0}$ of these   classical polynomials is orthogonal with respect to the measure $(1-t)^{\alpha} t^\beta dt,\ \alpha, \beta  > -1$ over the interval $[0,1].$ Therefore  the kernel $\Omega(x),\ x >0$ (1.11) is calculated accordingly, and we obtain

$$\Omega(x) = \int_{0}^1  (1-t)^{\alpha } t^{\beta-1} e^{-x/t} dt =  e^{-x} \int_{0}^\infty   t^{\alpha}\ (1+t)^{-\alpha- \beta-1} e^{-xt} dt= \Gamma(1+\alpha) e^{-x} U(1+\alpha, 1-\beta, x),\eqno(2.17)$$
where $U(a,b,z)$ is the Tricomi function \cite{nist}.  Hence we arrive at

{\bf Proposition 3.} {\it Let $\alpha > -1, \beta  > 0, \varphi, \psi$  be  nonnegative functions defined on  $\mathbb{R}_+$ which are related by the modified Laplace transform $(1.4)$.  Then under the condition 

$$\int_0^\infty  x^k e^{-x} \psi(x) dx < \infty, \quad k \in \mathbb{N}_0\eqno(2.18)$$
the sequence $\{P_n\}_{n\ge 0}$ of orthogonal polynomials with respect to the measure $e^{-x} U(1+\alpha, 1-\beta, x) \psi(x) dx$ over $\mathbb{R}_+$ is compositionally orthogonal in the sense of Jacobi relatively to the function $\varphi$, i.e.}

$$  \int_0^1  P_n(\theta) P_m(\theta) \{ \varphi(t)\}  (1-t)^\alpha  t^\beta dt= \Gamma(1+\alpha) \int_0^\infty  P_n(x) P_m(x)  e^{-x} U(1+\alpha, 1-\beta, x) \psi(x) dx = \delta_{n,m}.\eqno(2.19)$$

\begin{proof} Indeed, since

$$ P_n(\theta) P_m(\theta) \{ \varphi(t)\}  = P_n(\theta) P_m(\theta) \left\{ {1\over t} \int_0^\infty e^{-x/t} \psi(x) dx \right\} = {1\over t} \int_0^\infty e^{-x/t}  P_n(x) P_m(x) \psi(x) dx ,$$
where the consecutive differentiation under the integral sign is permitted owing to the estimate

$$ {1\over t} \int_0^\infty e^{-x/t} x^k \psi(x) dx  \le {1\over \delta}  \int_0^\infty e^{-x} x^k \psi(x) dx,\  0 < \delta \le t \le 1,\ k \in \mathbb{N}_0,$$
and the latter integral by $x$ is finite via (2.18).  Hence, taking into account (2.17), 

$$  \int_0^1  P_n(\theta) P_m(\theta) \{ \varphi(t)\}  (1-t)^\alpha  t^\beta dt=  \lim_{\delta \to 0+} \int_\delta^1  P_n(\theta) P_m(\theta) \{ \varphi(t)\}  (1-t)^\alpha  t^\beta dt$$

$$= \lim_{\delta \to 0+} \int_\delta^1  \int_0^\infty e^{-x/t}  P_n(x) P_m(x) \psi(x)  (1-t)^\alpha  t^{\beta-1} dx dt $$

$$=   \int_0^1  \int_0^\infty e^{-x/t}  P_n(x) P_m(x) \psi(x)  (1-t)^\alpha  t^{\beta-1} dx dt$$

$$  = \Gamma(1+\alpha) \int_0^\infty  P_n(x) P_m(x)  e^{-x} U(1+\alpha, 1-\beta, x) \psi(x) dx = \delta_{n,m},$$
where the interchange of the order of integration is possible by Fubini's theorem due to (2.18) and an elementary estimate of the Tricomi function

$$\Gamma(1+\alpha) U(1+\alpha, 1-\beta, x) = \int_{0}^\infty   t^{\alpha}\ (1+t)^{-\alpha- \beta-1} e^{-xt} dt \le   \int_{0}^\infty   t^{\alpha}\ (1+t)^{-\alpha- \beta-1} dt =
B(1+\alpha,\beta),$$
where $B(a,b)$ is the Euler beta function \cite{Bateman}, Vol. I.

\end{proof}

\section{Properties of Prudnikov's weight functions and their products}

In this section we will exhibit properties of the weight functions $\rho_{\nu+1}(x)\rho_\nu(x),\ \rho_\nu^2(x),\ x >0$, where $\rho_\nu(x)$ is defined by (2.2), (2.3), in order to involve them in the sequel to investigate the corresponding  orthogonal and multiple orthogonal polynomials.  In particular, we will establish their differential properties, integral representations and differential equations. Concerning the function $\rho_\nu$, we found in \cite{YAP}  the following integral representation in terms of Laguerre polynomials  

$${(-1)^n x^n\over n!}\  \rho_\nu(x)=   \int_0^\infty t^{\nu+n -1} e^{-t - x/t}  L_n^\nu(t) dt,\quad    n \in\mathbb{N}_{0}.\eqno(3.1)$$
It has a relationship with the  Riemann-Liouville fractional integral \cite{YaL}

$$ \left( I_{-}^\alpha  f \right) (x) \equiv \left( I_{-}^\alpha  f(x) \right)  = {1\over \Gamma(\alpha)} \int_x^\infty (t-x)^{\alpha-1} f(t) dt,\quad  {\rm Re} \alpha > 0,\eqno(3.2)$$
namely,  we get  the  formula 
$$\rho_\nu(x)= \left( I_{-}^\nu \rho_0 \right) (x),\ \nu >0.\eqno(3.3)$$
Further, the index law for fractional integrals immediately implies

$$ \rho_{\nu+\mu} (x)= \left( I_{-}^\nu \rho_\mu \right) (x)=   \left( I_{-}^\mu \rho_\nu \right) (x).\eqno(3.4)$$
The corresponding definition of the fractional derivative presumes the relation $ D^\mu_{-}= - D  I_{-}^{1-\mu}$.   Hence for the ordinary $n$-th derivative of $\rho_\nu$ we find

$$D^n \rho_\nu(x)= (-1)^n \rho_{\nu-n} (x),\quad n \in \mathbb{N}_0.\eqno(3.5)$$
 The  function $\rho_\nu$ possesses  the following recurrence relation (see \cite{YAP})

$$\rho_{\nu+1} (x) =    \nu \rho_\nu(x)+ x \rho_{\nu-1} (x),\quad \nu \in \mathbb{R}.\eqno(3.6)$$
In the operator form it can be written as follows

$$\rho_{\nu+1} (x) = \left( \nu - xD \right) \rho_\nu(x).\eqno(3.7)$$

Now we will derive the Mellin-Barnes representations for the product of functions $\rho_{\nu+1}\rho_\nu$ and for the square $\rho_\nu^2$. In fact, it can be done,  using the related formulas for the product of Macdonald functions.  Thus, appealing to Entries 8.4.23.31, 8.4.23.27 in \cite{PrudnikovMarichev}, Vol. III, we obtain

$$\rho_{\nu+1} (x) \rho_\nu (x) = {4^{-\nu} \over 2 i \sqrt\pi} \int_{\gamma-i\infty}^{\gamma +i\infty} \frac{\Gamma(   s+2\nu+1)  \Gamma(s+\nu) \Gamma(s)}
{\Gamma(s+\nu+1/2)}  (4x)^{-s} ds,\quad  x,\gamma > 0,\eqno(3.8) $$

$$\rho^2_\nu (x) = {4^{-\nu} \over i\sqrt\pi} \int_{\gamma-i\infty}^{\gamma +i\infty} \frac{\Gamma(   s+2\nu)  \Gamma(s+\nu) \Gamma(s)}
{\Gamma(s+\nu+1/2)}  (4x)^{-s} ds,\quad  x,\gamma > 0.\eqno(3.9) $$
Immediate consequences of these formulas are relationships of  the products of Prudnikov's weight functions $\rho_{\nu+1}\rho_\nu,\  \rho_\nu^2$ with $\rho_{2\nu+1},\ \rho_{2\nu}$, correspondingly. Indeed, it can be done via Entry 8.4.2.3   in \cite{PrudnikovMarichev}, Vol. III and the Parseval equality for the Mellin transform. Hence we deduce from (3.8), (3.9), respectively,

$$\rho_{\nu+1} (x) \rho_\nu (x)  = 4^{-\nu} \int_0^1 (1-t)^{-1/2} t^{\nu-1} \rho_{2\nu+1}\left({4x\over t}\right) dt$$

$$=  {x^{\nu} \sqrt\pi\over \Gamma(\nu+1/2)}  \int_0^1 (1-t)^{\nu-1/2} t^{-\nu-1} \rho_{\nu+1}\left({4x\over t}\right) dt,\eqno(3.10)$$

$$ \rho^2_\nu (x)  = 2^{1-2\nu} \int_0^1 (1-t)^{-1/2} t^{\nu-1} \rho_{2\nu}\left({4x\over t}\right) dt$$

$$= {2 x^\nu \sqrt\pi \over \Gamma(\nu+1/2)}  \int_0^1 (1-t)^{\nu-1/2} t^{-\nu-1} \rho_{\nu}\left({4x\over t}\right) dt .\eqno(3.11)$$
An ordinary differential equation for the function $\rho_{\nu+1}\rho_\nu$ is given by 

{\bf Proposition 4.} {\it The function $u_\nu = \rho_{\nu+1}\rho_\nu$ satisfies the following third order differential equation}

$$x^2{d^3u_\nu\over dx^3} + x(2- 3\nu) {d^2 u_\nu\over dx^2} + 2 (\nu( \nu-1) -2x) {d u_\nu\over dx} + 2(2\nu-1) u_\nu(x) = 0.\eqno(3.12)$$

\begin{proof}  Differentiating $u_\nu$ with the use of (3.5), we have

$$ {d u_\nu\over dx} = - \rho_\nu^2(x) -   \rho_{\nu+1}(x) \rho_{\nu-1} (x).\eqno(3.13)$$
Multiplying both sides of (3.13) by $x$ and employing the recurrence relation (3.6), we get

$$x {d u_\nu\over dx} = - x \rho_\nu^2(x) + \nu  \rho_{\nu+1}(x) \rho_{\nu} (x) -  \rho^2_{\nu+1}(x).$$
Differentiating both sides of the latter equality, using (3.6), (3.13) and the notation $u_\nu= \rho_{\nu+1}\rho_\nu$, we obtain

$${d\over dx} \left(x {d u_\nu\over dx}\right)  = -  (1+2\nu) \rho_\nu^2(x)  +\nu {d u_\nu\over dx}+ 4u_\nu (x).\eqno(3.14)$$
Differentiating (3.14) and multiplying the result by $x$, we recall (3.6) to find

$$x  {d^2\over dx^2} \left(x {d u_\nu\over dx}\right)  =  - 2\nu  (1+2\nu) \rho^2_\nu(x)   +\nu x {d^2 u_\nu\over dx^2}+ 4x {d u_\nu\over dx} + 2(1+2\nu) u_\nu(x).\eqno(3.15)$$
Hence, expressing $ (1+2\nu) \rho^2_\nu(x)$ from (3.14) and fulfilling the differentiation on the left-hand side of (3.15), we get (3.12).

\end{proof}

Concerning the differential equation for the function $\rho_\nu^2$, we prove

{\bf Proposition 5.} {\it The function $h_\nu = \rho^2_\nu$ satisfies the following third order differential equation}

$$x^2{d^3h_\nu\over dx^3} +  3x(1- \nu) {d^2 h_\nu\over dx^2} +  (2\nu^2+1- 3\nu- 4x) {d h_\nu\over dx} + 2(2\nu-1) h_\nu(x) = 0.\eqno(3.16)$$

\begin{proof}  Since $x \left(\rho^2_\nu(x) \right)^\prime = -2x \rho_{\nu-1}(x) \rho_{\nu} (x) = 2\nu \rho^2_\nu(x) - 2 \rho_{\nu+1}(x) \rho_{\nu} (x) = 2\nu h_\nu(x)  - 2 \rho_{\nu+1}(x) \rho_{\nu} (x)  $  we derive, owing to (3.5), (3.6)

$$x {d\over dx} \left(  x {d h_\nu\over dx} - 2\nu h_\nu(x) \right) = -2 x {d\over dx} \left( \rho_{\nu+1}(x) \rho_{\nu} (x)
\right)$$

$$ = 2x\ h_\nu(x) - 2\nu   \rho_{\nu+1}(x) \rho_{\nu} (x) + 2 \rho_{\nu+1}^2(x)=  2 \left( x - \nu^2 \right)  h_\nu(x)   +\nu x  {d h_\nu\over dx}  + 2 \rho_{\nu+1}^2(x).$$
Hence one more differentiation yields

$$  {d\over dx} \left( x {d\over dx} \left(  x {d h_\nu\over dx} - 2\nu h_\nu(x) \right) \right) =    {d\over dx} \left(  2 \left( x - \nu^2 \right)  h_\nu(x)   +\nu x  {d h_\nu\over dx}\right) - 4  \rho_{\nu+1}(x) \rho_{\nu} (x).$$
But from the beginning of the proof we find 

$$2 \rho_{\nu+1}(x) \rho_{\nu} (x) = 2\nu  h_\nu(x) -  x {d h_\nu\over dx}.\eqno(3.17)$$ 
Therefore, substituting this expression into the previous equality, we fulfil the differentiation to arrive at (3.16) and to complete the proof.    

\end{proof}

{\bf Corollary 1}.  {\it The following recurrence relations between functions $u_\nu,\ h_\nu$ hold}

$$ u_\nu=  \nu h_\nu + x u_{\nu-1},\eqno(3.18)$$ 

$$h_{\nu+1} = \nu^2 h_\nu + 2 x \nu u_{\nu-1} + x^2 h_{\nu-1} = 2\nu u_\nu+ x^2 h_{\nu-1} -  \nu^2 h_\nu .\eqno(3.19)$$

\begin{proof} Equality (3.18) is a direct consequence of (3.17) and (3.5).  Equalities (3.19), in turn,  are obtained, taking squares of both sides of (3.6) and employing (3.18).

\end{proof}

\section{Orthogonal polynomials with  $\rho^2_\nu$ weight function}

 The object of this section is  to  characterize the sequence of orthogonal polynomials $\{P_n\}_{n\ge 0}$, satisfying the orthogonality conditions

$$\int_{0}^{\infty}  P_{m} (x) P_{n}(x)  \rho^2_\nu(x) dx = \delta_{n,m}, \quad \nu > - {1\over 2}.\eqno(4.1)$$
Clearly,  up to a normalization factor conditions (4.1) are equivalent to the equalities 

$$\int_{0}^{\infty}  P_{n} (x)  \rho^2_\nu(x)   x^m dx = 0,\quad m=0,1,\dots, n-1,\ \quad n \in \mathbb{N}.\eqno(4.2)$$
The moments of the weight $\rho^2_\nu(x)$ can be obtained immediately from the Mellin-Barnes representation (3.9), treating it as the inverse Mellin transform \cite{Tit}. Hence we get 

$$ \int_{0}^{\infty}   \rho^2_\nu(x)   x^\mu dx =  \sqrt\pi \ \frac{\Gamma(   1+\mu+2\nu)  \Gamma(1+\mu+\nu) \Gamma(1+\mu)}
{ 2^{1+2(\nu+\mu)}\  \Gamma(\mu+\nu+3/2)}, \  \mu > \max\{-1, -1-\nu,\ 1-2\nu\}.\eqno(4.3)$$
Furthermore, the sequence $ \left\{P_n\right\}_{n \ge 0}$ satisfies the 3-term recurrence relation in the form

$$x  P_n (x) = A_{n+1}  P_{n+1} (x) + B_n  P_n (x) + A_n  P_{n-1} (x),\eqno(4.4)$$
where $P_{-1}^{\nu} (x)\equiv 0,\  P_n(x)= \sum_{k=0}^n a_{n,k} x^k,\  a_{n,n} \neq 0$ and

$$A_{n+1}= {a_n\over a_{n+1} }, \quad B_{n}= {b_n\over a_n} -  {b_{n+1}\over a_{n+1}},\quad  a_n\equiv a_{n,n},\  b_n\equiv a_{n,n-1}.$$
Basing on properties of the functions $\rho^2_\nu,\ \rho_{\nu+1}\rho_\nu$ and orthonormality (4.1), we establish 

{\bf Proposition 6}.  {\it Let $\nu > -1$. The following formulas hold}

$$ \int_{0}^{\infty}   P^2_{n}(x)  \rho_{\nu+1}(x)\rho_\nu(x)  dx =  {1\over 2} + \nu + n,\eqno(4.5)$$

$$ \int_{0}^{\infty}   P^2_{n}(x) \ x \rho_{\nu}(x) \rho_{\nu-1}(x) dx =   {1\over 2}  + n,\eqno(4.6)$$

$$ \int_{0}^{\infty}   P^2_{n}(x) \ x^2 \rho_{\nu}(x) \rho_{\nu-2}(x) dx =  B_n- (\nu-1) \left( {1\over 2}  + n\right),\eqno(4.7)$$

$$ \int_{0}^{\infty}   P^2_{n}(x)  \rho_{\nu+2}(x)\rho_\nu(x)  dx = B_n+ (\nu+1)\left( {1\over 2} + \nu + n\right).\eqno(4.8)$$

\begin{proof} In fact, recurrence relation (3.6) and orthonormality (4.1) implies

$$ \int_{0}^{\infty}   P^2_{n}(x)  \rho_{\nu+1}(x)\rho_\nu(x)  dx  = \nu +  \int_{0}^{\infty}   P^2_{n}(x)\  x \rho_{\nu}(x)\rho_{\nu-1}(x)  dx.\eqno(4.9)$$
 Hence, integrating by parts in the latter integral and eliminating integrated terms by virtue of the asymptotic behavior  (2.4), (2.5) , we find
 
 $$\int_{0}^{\infty}   P^2_{n}(x) \ x \rho_{\nu}(x)\rho_{\nu-1}(x)  dx = \int_0^\infty  \left( {1\over 2} P^2_{n}(x) +  x P_{n}(x) P_n^\prime (x) \right) \rho^2_{\nu}(x) dx.\eqno(4.10) $$
Therefore, we have from (4.9), (4.10)  and orthogonality (4.1)

$$ \int_{0}^{\infty}   P^2_{n}(x)  \rho_{\nu+1}(x)\rho_\nu(x)  dx = {1\over 2} + \nu + n,$$
which proves (4.5).   Concerning equality (4.6), it is a direct consequence of (3.6), (4.1).  The same idea is to prove (4.7), (4.8), employing (4.3) as well.

\end{proof}

The composition orthogonality which is associated with the sequence  $\{P_n\}_{n\ge 0}$ is given by

{\bf Theorem 1.} {\it Let $\nu > -1/2$. The sequence $\{P_n\}_{n\ge 0}$ is compositionally orthogonal in the sense of Laguerre relatively to the function $t^\nu e^t \Gamma(-\nu,t)$, where $\Gamma(\mu,z)$ is the incomplete gamma function, i.e.}

$$\ \int_{0}^{\infty}   t^{\nu }   e^{-t}  P_n\left(\theta\right)  P_m\left(\theta\right)\left(  t^\nu e^t \Gamma(-\nu, t) \right) dt = { \delta_{n,m}\over \Gamma(1+\nu)}.\eqno(4.11)$$

\begin{proof} Writing $P_n$ explicitly, we appeal to the integral representation (3.1) for $\rho_\nu(x)$ to write the left-hand side of (4.2) in the form

$$ \int_{0}^{\infty}  P_{n} (x)  \rho^2_\nu(x)   x^m dx =  (-1)^m m! \sum_{k=0}^n a_{n,k}  (-1)^k k!  \int_{0}^{\infty} \int_{0}^{\infty}   t^{\nu+k -1} e^{-t - x/t}  L_k^\nu(t) dt$$

$$\times \int_{0}^{\infty}   y^{\nu+m -1} e^{-y - x/y}  L_m^\nu(y) dy dx.\eqno(4.12) $$
The Fubini theorem allows to interchange the order of integration in (4.12) due to the convergence of the following integral for all $r_1, r_2  \in \mathbb{N}_0$

$$\int_{0}^{\infty} \int_{0}^{\infty}   t^{\nu+k+r_1 -1} e^{-t - x/t} dt  \int_{0}^{\infty}   y^{\nu+m +r_2 -1} e^{-y - x/y}  dy\ dx $$

$$=   \int_{0}^{\infty}   \int_{0}^{\infty}   t^{\nu+k+r_1 }   y^{\nu+m +r_2} e^{-t-y } {dt dy\over t+ y} \le {1\over 2} \int_{0}^{\infty}   \int_{0}^{\infty}   t^{\nu+k+r_1 -1/2 }   y^{\nu+m +r_2-1/2 } e^{-t-y } dt dy$$

$$ =   {1\over 2} \Gamma(\nu+k+r_1 +1/2)  \Gamma(\nu+m+r_2 +1/2).$$
Therefore after integration with respect to $x$ we find from (4.12)

$$ \int_{0}^{\infty}  P_{n} (x)  \rho^2_\nu(x)   x^m dx =  (-1)^m m! \sum_{k=0}^n a_{n,k}  (-1)^k k!  \int_{0}^{\infty} \int_{0}^{\infty}   t^{\nu+k}   y^{\nu+m } e^{-t-y}  L_k^\nu(t)  L_m^\nu(y) {dt dy \over t+y}.\eqno(4.13) $$
However the Rodrigues formula for Laguerre polynomials and Viskov-type identity  (1.3) suggest to write the right-hand side of (4.13) as follows

$$  (-1)^m m! \sum_{k=0}^n a_{n,k}  (-1)^k k!  \int_{0}^{\infty} \int_{0}^{\infty}   t^{\nu+k}   y^{\nu+m } e^{-t-y}  L_k^\nu(t)  L_m^\nu(y) {dt dy \over t+y}$$

$$=   (-1)^m \sum_{k=0}^n a_{n,k}  (-1)^k k!  \int_{0}^{\infty}   t^{\nu+k}   e^{-t}   L_k^\nu(t)  \int_{0}^{\infty}   \theta^m \left( y^\nu e^{-y} \right) {dy dt \over t+y} . \eqno(4.14)$$
The inner integral with respect to $y$ on the right-hand side of the latter equality can be treated via integration by parts, and we obtain

$$ \int_{0}^{\infty}   \theta^m \left( y^\nu e^{-y} \right) {dy  \over t+y}  = (-1)^m \int_{0}^{\infty}    y^\nu e^{-y} \theta^m \left( {1\over t+y} \right) dy.\eqno(4.15)$$
Working out the differentiation, we find

$$ \theta^m \left( {1\over t+y} \right) =  y^m {d^m\over dy^m} \left( \sum_{k=0}^m \binom{m}{k} (-1)^{m-k} t^{m-k} (t+y)^{k-1}\right)  $$

$$=  (-1)^m t^m y^m {d^m\over dy^m} \left(  {1\over t+y}\right) = {m! \  t^m y^m \over (t+y)^{m+1} }.$$
Consequently, combining with (4.14), (4.15), equality (4.13) becomes

$$ \int_{0}^{\infty}  P_{n} (x)  \rho^2_\nu(x)   x^m dx  =    m! \sum_{k=0}^n a_{n,k}  (-1)^k k!  \int_{0}^{\infty}   t^{\nu+k+m }   e^{-t}   L_k^\nu(t)  \int_{0}^{\infty} 
  {y^{\nu+m}  e^{-y}  \over (t+y)^{m+1} } dy dt.\eqno(4.16) $$
Meanwhile, the integral with respect to $y$ in (4.16) is calculated in  (2.17) in terms of the Tricomi function up to a simple change of variables.  Thus we get

 $$ \int_{0}^{\infty}  P_{n} (x)  \rho^2_\nu(x)   x^m dx  =    m!  \Gamma(\nu+m+1) \sum_{k=0}^n a_{n,k}  (-1)^k k!  \int_{0}^{\infty}   t^{2\nu+k+m }   e^{-t}   L_k^\nu(t)  U(\nu+m+1, 1+\nu, t) dt.\eqno(4.17) $$
But appealing to the differential formula 13.3.24 in \cite{nist} for the Tricomi function, we have

$$  m! \ \Gamma(\nu+m+1)  t^{\nu+m}  U(\nu+m+1, 1+\nu, t) = \Gamma(\nu+1) \theta^m\left( t^\nu U(1+\nu,1+\nu,t)\right)$$

$$ = \Gamma(\nu+1) \theta^m\left( t^\nu e^t \Gamma(-\nu, t) \right),\eqno(4.18)$$
where $\Gamma(a,z)$ is the incomplete gamma function

$$ \Gamma(a,z) = \int_z^\infty e^{-u} u^{a-1} du.\eqno(4.19)$$
Hence,  plugging this result in (4.17), recalling the Rodrigues  formula for Laguerre polynomials and integrating by parts on its right-hand side, we  combine  with (4.2) to  end up with equalities   

$$ \int_{0}^{\infty}  P_{n} (x)  \rho^2_\nu(x)   x^m dx = \Gamma(\nu+1) $$

$$\times  \int_{0}^{\infty}   t^{\nu }   e^{-t}  P_n\left(\theta\right) \theta^m \left(  t^\nu e^t \Gamma(-\nu, t) \right) dt = 0,\quad m=0,1,\dots, n-1,\  n \in \mathbb{N},\eqno(4.20)$$
which yield (4.11) and complete  the proof of Theorem 1. 

\end{proof}

Further, the right-hand side of the first equality in (4.20) can be rewritten via integration by parts as follows

$$ \int_{0}^{\infty}   t^{\nu+m }   e^{-t}    L_m^\nu(t) P_n\left(\theta\right)  \left(  t^\nu e^t \Gamma(-\nu, t) \right) dt = 0,\quad m=0,1,\dots, n-1,\  n \in \mathbb{N}.\eqno(4.21)$$
Then we expand the function $ F_n(t)= \theta^n \left(  t^\nu e^t \Gamma(-\nu, t) \right)$  in a series of Laguerre polynomials 

$$ F_n(t)=   \sum_{r=0}^\infty c_{n,r} L_r^\nu(t),\eqno(4.22)$$
where   the coefficients $c_{n,r}$ are calculated in terms of the ${}_3F_2$-hypergeometric functions at  unity with the aid of Entry 3.31.12.1 in \cite{Brychkov}. Precisely, we get via (4.18)

$$c_{n,r}= {r!\over \Gamma(1+\nu+r)} \int_0^\infty t^\nu e^{-t}  L_r^\nu(t) \theta^n \left(  t^\nu e^t \Gamma(-\nu, t) \right) dt$$

$$= {r! n! (1+\nu)_n \over \Gamma(1+\nu+r)} \int_0^\infty t^{2\nu+n}  e^{-t}  L_r^\nu(t) U(\nu+n+1, 1+\nu, t)  dt$$

$$=    {(1+\nu)_n\over \Gamma(1+\nu)}   \left[ {n! r! \ \Gamma(\nu) \over (1+\nu)_r}  \  {}_3F_2 \left(  n+1, r+1, \nu+n+1; 1-\nu,  1; \ 1 \right) \right.$$

$$\left.  +   \Gamma(2\nu+n+1) \Gamma(-\nu)\  {}_3F_2 \left( \nu+ r+1,  \nu+n+1,  2\nu+n+1; 1+\nu,  1+\nu; \ 1 \right)\right],\eqno(4.23)$$
where $(a)_z$ is the Pochhammer symbol \cite{Bateman}, Vol. I and it is valid for nonnegative integers $\nu$ by continuity.  Hence after substitution the series (4.22) into (4.21) these orthogonality conditions  take the form

$$  \int_{0}^{\infty}   t^{\nu+m }   e^{-t}    L_m^\nu(t) \sum_{k=0}^n a_{n,k} \sum_{r=0}^{2m}  c_{k,r} L_r^\nu(t)   dt = 0,\quad m=0,1,\dots, n-1,\  n \in \mathbb{N}.\eqno(4.24)$$
Calculating the integral in (4.24) via relation (2.19.14.8) in \cite{PrudnikovMarichev}, Vol. II we obtain

 $$d_{m,r} =  \int_{0}^{\infty}   t^{\nu+m }   e^{-t}    L_m^\nu(t)  L_r^\nu(t)   dt = { (-1)^r   \over  r! } \  (1+\nu)_r\  \Gamma(1+\nu+m) $$
 
 $$\times  {}_3F_2 \left( -r, \nu+ m+1,  m+1; 1+\nu,  1; \ 1 \right).\eqno(4.25)$$ 
Hence, equalities (4.24) become the linear system of $n$ algebraic equations with $n+1$ unknowns  

$$\sum_{k=0}^n a_{n,k} f_{k,m} = 0,\quad m=0,1,\dots, n-1,\  n \in \mathbb{N},\eqno(4.26)$$
where

$$f_{k,m} = \sum_{r=0}^{2m}  c_{k,r} d_{m,r}.\eqno(4.27)$$ 
Consequently,  explicit values of the coefficients $a_{n,k},\ k =1,2,\dots, n$ can be expressed via Cramer's rule in terms of the free coefficient $a_{n,0}$ as follows  

$$  a_{n,k} = -    a_{n,0}\  { D_{n,k}  \over D_{n}},\quad   k = 1,\dots, n,\eqno(4.28)$$
where

$$ D_{n}= \begin{vmatrix}  

f_{1, 0}  & f_{2, 0}&  \dots&  \dots&   f_{n, 0} \\

f_{1, 1} &  \dots &   \dots&    \dots&    f_{n, 1} \\
 
\dots&  \dots &   \dots&   \dots&   \dots \\
 
 \vdots  &   \ddots  &  \ddots & \ddots&   \vdots\\
  
 f_{1, n-1} &   \dots&  \dots&   \dots&    f_{n, n-1}\\

 \end{vmatrix},\eqno(4.29)$$

$$ D_{n,k}= \begin{vmatrix}  

f_{1, 0}  & \dots & f_{k-1, 0}&  f_{0,0} &   f_{k+1, 0} & \dots&   f_{n, 0}  \\

f_{1, 1} &  \dots & f_{k-1, 1}&  f_{0,1} &   f_{k+1, 1} &  \dots&   f_{n, 1}  \\

 \dots&  \dots &   \dots&   \dots&  \dots&   \dots&  \dots \\
 
 \vdots  &   \ddots  &  \ddots & \vdots&  \ddots  &  \ddots & \vdots\\
  
 \vdots  &   \ddots  &  \ddots & \vdots&  \ddots  &  \ddots & \vdots\\

 \vdots  &   \ddots  &  \ddots & \vdots&  \ddots  &  \ddots & \vdots\\

f_{1, n-1} &  \dots & f_{k-1, n-1}&  f_{0,n-1} &   f_{k+1, 1} &  \dots&   f_{n, n-1}  \\

 \end{vmatrix}.\eqno(4.30)$$
The free  coefficient can be determined, in turn,  from the orthogonality conditions (4.1), (4.2), which imply the formula

$$\int_{0}^{\infty}  P_{n} (x)   \rho^2_\nu(x) \  x^n dx = {1\over a_{n,n}}.\eqno(4.31)$$ 
Therefore from (4.28) and (4.3)  we derive

$$ {1\over a_{n,n}} =  -   { a_{n,0} \sqrt\pi \over D_n  }\  \frac{\Gamma(   1+2\nu)  \Gamma(1+\nu) }
{ 2^{1+2\nu}\  \Gamma(\nu+3/2)} \ \sum_{k=0}^n  D_{n,k} \ (n+k)! \ \frac{(   1+2\nu)_{n+k}   (1+\nu)_{n+k} }
{ 4^{n+k}\  (\nu+3/2)_{n+k}},\eqno(4.32)$$
where $D_{n,0}\equiv - D_n.$  Hence

$$  a_{n,0}  =    \pm \  { D_n  \over [ D_{n,n} ]^{1/2} }\left[ \sqrt\pi \ \frac{\Gamma(   1+2\nu)  \Gamma(1+\nu) }
{ 2^{1+2\nu}\  \Gamma(\nu+3/2)} \ \sum_{k=0}^n  D_{n,k} \ (n+k)! \ \frac{(   1+2\nu)_{n+k}   (1+\nu)_{n+k} }
{ 4^{n+k}\  (\nu+3/2)_{n+k}} \right]^{-1/2},\eqno(4.33)$$
where the sign can be  chosen  accordingly, making positive expressions under the square roots.   Assuming also the positivity of the leading coefficient $a_{n,n}$ we have its value, correspondingly,

$$   a_{n,n}  =    \mp \   [ D_{n,n}]^{1/2} \left[  \sqrt\pi \ \frac{\Gamma(   1+2\nu)  \Gamma(1+\nu) }
{ 2^{1+2\nu}\  \Gamma(\nu+3/2)} \ \sum_{k=0}^n  D_{n,k} \ (n+k)! \ \frac{(   1+2\nu)_{n+k}   (1+\nu)_{n+k} }
{ 4^{n+k}\  (\nu+3/2)_{n+k}}  \right]^{-1/2}.\eqno(4.34)$$

{\bf Theorem 2}.  {\it Let $\nu > -1/2$. The sequence of orthogonal polynomials $\left\{P_n\right\}_{n \ge 0}$  can be expressed explicitly, where the coefficients $a_{n,k},\ k=1,2,\dots, n$ are calculated by relations  $(4.28)$ and the free term $a_{n,0}$ is defined by the equality $(4.33)$. Besides,  it satisfies the  3-term recurrence relation $(4.4)$, where }

$$A_{n+1}= { a_{n,0} \  D_{n+1} \ D_{n,n}   \over a_{n+1,0}\  D_n \ D_{n+1,n+1}},  \quad  B_n= { D_{n,n-1}  \over D_{n,n}} -  { D_{n+1,n}  \over D_{n+1,n+1}}.\eqno(4.35)$$
An analog of the Rodrigues formula for the sequence $ \left\{P_n\right\}_{n \in \mathbb{N}_0}$ can be established, appealing to the   integral representation  of an arbitrary polynomial in terms of the associated polynomial of degree $2n$ (see \cite{YAP}). Hence, employing (2.2), we deduce

$$P_n(x)= {1\over \rho_\nu(x)}   \int_0^\infty t^{\nu-1} e^{-t -x/t  }  q_{2n}(t) dt = {1\over \rho^2_\nu(x)}     \int_0^\infty u^{\nu-1} e^{-u -x/u  }  du \int_0^\infty t^{\nu-1} e^{-t -x/t  }  q_{2n}(t) dt $$

$$ = {(-1)^n \over \rho^2_\nu(x)}    {d^n\over dx^n } \int_0^\infty  \int_0^\infty (u t)^{\nu+n-1} e^{-u -t -x(1/u +1/t)  }  q_{2n}(t) {dt du\over (u+t)^n},\eqno(4.36)$$
where
$$ q_{2n}(x) =  \sum_{k=0}^{n} a_{n,k} (-1)^k k! x^k L_k^\nu(x)\eqno(4.37)$$
and the $k$-th differentiation under the integral sign is permitted due to the estimate

$$\int_0^\infty  \int_0^\infty (u t)^{\nu+k-1} e^{-u -t -x(1/u +1/t)  } \left| q_{2n}(t)\right|  {dt du\over (u+t)^k}$$

$$ \le 2^{-k}  \int_0^\infty   u^{\nu+ k/2 -1}  e^{-u} du  \int_0^\infty   t^{\nu+ k/2 -1}  e^{-t} \left| q_{2n}(t)\right|  dt $$

$$=   2^{-k}  \Gamma\left( \nu+ {k\over 2} \right)  \int_0^\infty   t^{\nu+ k/2 -1}  e^{-t} \left| q_{2n}(t)\right|  dt   < \infty,\quad  k = 1,\dots, n.$$
Then writing

$${(ut)^n\over (u+t)^n} = {1\over (n-1)!} \int_0^\infty e^{ - (1/u+1/t) y} y^{n-1} dy,$$
we get from (4.36)

$$P_n(x)=  {(-1)^n \over (n-1)!  \ \rho^2_\nu(x)}    {d^n\over dx^n } \int_0^\infty  \int_0^\infty   \int_0^\infty  y^{n-1} (u t)^{\nu-1} e^{-u -t -(x+y)/u -(x+y)/t  }  q_{2n}(t) dt du dy$$

$$=  {(-1)^n \over (n-1)!  \ \rho^2_\nu(x)}    {d^n\over dx^n }  \int_0^\infty   \int_0^\infty  y^{n-1} t^{\nu-1} e^{ -t  -(x+y)/t  }  q_{2n}(t) \rho_{\nu}(x+y) dt dy.\eqno(4.38)$$
Then, expressing $q_{2n}$ in terms of the  Laguerre polynomials 

$$  q_{2n}(x)= \sum_{k=0}^{2n} h_{2n,k} L_k^\nu(x),\eqno(4.39)$$
where the coefficients $h_{2n,k}$ are calculated by virtue of (4.28),  (4.37) and relation (2.19.14.15) in \cite{PrudnikovMarichev}, Vol. II, namely, 

$$ h_{2n,k} =  {k!\over \Gamma(1+\nu+k)} \int_0^\infty t^\nu e^{-t} L_k^\nu(t)  q_{2n}(t) dt $$

$$= -  {k! \ a_{n,0} \over  D_n\ \Gamma(1+\nu+k)} \sum_{r=0}^n D_{n,r} (-1)^r r!  \int_0^\infty t^{\nu+r}  e^{-t} L_k^\nu(t)  L_r^\nu(t) dt $$

$$= -  { a_{n,0} \over  D_n} \sum_{r=0}^n D_{n,r} \  r! \ (1+\nu)_r  \ {}_3F_2 \left(-k,\ 1+\nu+r,\ 1+r;\ 1+\nu,\ 1;\ 1 \right). $$
So, 

$$ h_{2n,k} = -  { a_{n,0} \over  D_n} \sum_{r=0}^n D_{n,r} \  r! \ (1+\nu)_r  \ {}_3F_2 \left(-k,\ 1+\nu+r,\ 1+r;\ 1+\nu,\ 1;\ 1 \right),\eqno(4.40) $$
and  the values of the generalized hypergeometric function can be simplified via relations (7.4.4; 90,91,92,93) in \cite{PrudnikovMarichev}, Vol. III. Precisely, we get for $k =0,1\dots, n$

$$ {}_3F_2 \left(-k,\ 1+\nu+r,\ 1+r;\ 1+\nu,\ 1;\ 1 \right) = 0,\quad k > 2r,\eqno(4.41)$$

$$ (1+\nu)_r\  {}_3F_2 \left(-k,\ 1+\nu+r,\ 1+r;\ 1+\nu,\ 1;\ 1 \right) = {  (2r) ! \over  r!  },\quad k=2r,$$ 

$$ (1+\nu)_r\  {}_3F_2 \left(-k,\ 1+\nu+r,\ 1+r;\ 1+\nu,\ 1;\ 1 \right) = -  {  (2r-1) ! (\nu+ 3r ) \over  (r-1)! },  \quad k=2r-1, $$

$$ (1+\nu)_r\  {}_3F_2 \left(-k,\ 1+\nu+r,\ 1+r;\ 1+\nu,\ 1;\ 1 \right) = {  (2(r-1)) ! \over 2\  r! }
\left( 2 r^2 (2r+\nu-1)(2r-1) \right.$$

$$\left. +  r(r-1)(r+\nu-1)(r+\nu ) \right),\quad k=2(r-1).$$ 
Therefore, returning to (4.38) and minding (4.39), (4.40), (4.41), (3.1), (3.2) we find

$$P_n(x)= {(-1)^{n+1}  a_{n,0}  \over  D_n\  (n-1)!  \ \rho^2_\nu(x)}    {d^n\over dx^n }   \sum_{k=0}^{2n}   \sum_{r=0}^n   D_{n,r} \  r! \ (1+\nu)_r  \ {}_3F_2 \left(-k,\ 1+\nu+r,\ 1+r;\ 1+\nu,\ 1;\ 1 \right) $$

$$\times \int_0^\infty   \int_0^\infty  y^{n-1} t^{\nu-1} e^{ -t  -(x+y)/t  } L_k^\nu(t)  \rho_{\nu}(x+y) dt dy$$

$$= {(-1)^{n+1}  a_{n,0}  \over  D_n\  (n-1)!  \ \rho^2_\nu(x)}    {d^n\over dx^n }   \sum_{k=0}^{2n} {1\over k!}   \sum_{r=0}^n   D_{n,r} \  r! \ (1+\nu)_r  \ {}_3F_2 \left(-k,\ 1+\nu+r,\ 1+r;\ 1+\nu,\ 1;\ 1 \right) $$

$$\times \int_0^\infty    y^{n-1}  {d^k\over dx^k} \left(  (x+y)^k \rho_\nu(x+y) \right) \rho_{\nu}(x+y) dy$$

$$= {(-1)^{n+1}  a_{n,0}  \over  D_n\  (n-1)!  \ \rho^2_\nu(x)}    {d^n\over dx^n }   \sum_{k=0}^{2n} {1\over k!}  \sum_{r=0}^n   D_{n,r} \  r! \ (1+\nu)_r  \ {}_3F_2 \left(-k,\ 1+\nu+r,\ 1+r;\ 1+\nu,\ 1;\ 1 \right) $$

$$\times  \int_x^\infty  (y-x)^{n-1}  {d^k\over dy^k} \left(y^k \rho_{\nu}(y)\right) \rho_\nu(y)  dy =  - { a_{n,0}  \over  D_n \ \rho_\nu(x)}    \sum_{k=0}^{2n}  {d^k\over dx^k}  \left( x^k \rho_{\nu}(x)\right)  \sum_{r=0}^n   D_{n,r} \  r! $$

$$\times \ {(1+\nu)_r\over k!}   \ {}_3F_2 \left(-k,\ 1+\nu+r,\ 1+r;\ 1+\nu,\ 1;\ 1 \right)$$

$$=  - { a_{n,0}  \over  D_n \ \rho_\nu(x)}   \sum_{r=0}^n    D_{n,r} \  r!    \ (1+\nu)_r \sum_{k=0}^{2r}    {d^k\over dx^k}  \left( x^k \rho_{\nu}(x)\right)  {1\over k!} \ {}_3F_2 \left(-k,\ 1+\nu+r,\ 1+r;\ 1+\nu,\ 1;\ 1 \right).$$

Thus it proves 

{\bf Theorem 3}. {\it Let $\nu > -1/2,\ n \in \mathbb{N}_0$.  Orthogonal polynomials $P_n$ satisfy  the  Rodrigues-type  formula 

$$P_n(x)= - { a_{n,0}  \over  D_n \ \rho_\nu(x)}   \sum_{r=0}^n    D_{n,r}   r! (1+\nu)_r \sum_{k=0}^{2r} {1\over k!}  {d^k\over dx^k}  \left( x^k \rho_{\nu}(x)\right){}_3F_2 \left(-k,\ 1+\nu+r,\ 1+r;\ 1+\nu,\ 1;\ 1 \right),\eqno(4.42)$$
where  $a_{n,0}$ is defined by  $(4.33)$ and $D_n,\ D_{n,r}$ by $(4.29), (4.30)$, respectively.} 

{\bf Corollary 2}. {\it  Orthogonal polynomials  $P_n$  have the following representation

$$P_n(x)= - { a_{n,0}  \over  D_n \ }   \sum_{r=0}^n    D_{n,r}   r! (1+\nu)_r\left[  \sum_{k=0}^{r} {A_{k,k-1} (x) \over (2k)!}   {}_3F_2 \left(-2k,\ 1+\nu+r,\ 1+r;\ 1+\nu,\ 1;\ 1 \right)\right.$$

$$\left. +  \sum_{k=0}^{r-1} {A_{k,k} (x) \over (2k+1)!}   {}_3F_2 \left(-2k-1,\ 1+\nu+r,\ 1+r;\ 1+\nu,\ 1;\ 1 \right) \right], \eqno(4.43)$$
where $ A_{k,k-1},  A_{k,k}$ are the type $1$ multiple orthogonal polynomials of degree $k$,  associated with the vector of weight functions $(\rho_\nu, \rho_{\nu+1})$ }.

\begin{proof} In fact, we write (4.42) in the form

 $$P_n(x)= - { a_{n,0}  \over  D_n \ \rho_\nu(x)}   \sum_{r=0}^n    D_{n,r}   r! (1+\nu)_r \left[ \sum_{k=0}^{r} {1\over (2k)!}  {d^{2k} \over dx^{2k}}  \left( x^{2k} \rho_{\nu}(x)\right){}_3F_2 \left(-2k,\ 1+\nu+r,\ 1+r;\ 1+\nu,\ 1;\ 1 \right)\right.$$

$$\left. +  \sum_{k=0}^{r-1} {1\over (2k+1)!}  {d^{2k+1} \over dx^{2k+1}}  \left( x^{2k+1} \rho_{\nu}(x)\right){}_3F_2 \left(-2k-1,\ 1+\nu+r,\ 1+r;\ 1+\nu,\ 1;\ 1 \right)\right].\eqno(4.44)$$
Meanwhile,  appealing to the Rodrigues formulas for the type $1$  multiple orthogonal polynomials  associated with the vector of weight functions $(\rho_\nu, \rho_{\nu+1})$ (see in \cite{AsscheYakubov2000}), it gives 

$$ {d^{2k} \over dx^{2k}}  \left( x^{2k} \rho_{\nu}(x)\right) = A_{k,k-1} (x)\rho_\nu(x) + B_{k,k-1} (x)\rho_{\nu+1}(x),$$

$$ {d^{2k+1} \over dx^{2k+1}}  \left( x^{2k+1} \rho_{\nu}(x)\right) = A_{k,k} (x)\rho_\nu(x) + B_{k,k} (x)\rho_{\nu+1}(x),$$
where $  A_{k,k-1},  B_{k,k-1}$ are polynomials of degree $k, k-1$, respectively, and  $  A_{k,k},  B_{k,k}$ are polynomials of degree $k$.   Therefore, substituting these expressions into (4.44), we obtain

$$P_n(x)= - { a_{n,0}  \over  D_n \ }   \sum_{r=0}^n    D_{n,r}   r! (1+\nu)_r\left[  \sum_{k=0}^{r} {A_{k,k-1} (x) \over (2k)!}   {}_3F_2 \left(-2k,\ 1+\nu+r,\ 1+r;\ 1+\nu,\ 1;\ 1 \right)\right.$$

$$\left. +  \sum_{k=0}^{r-1} {A_{k,k} (x) \over (2k+1)!}   {}_3F_2 \left(-2k-1,\ 1+\nu+r,\ 1+r;\ 1+\nu,\ 1;\ 1 \right) \right] $$

$$  - { a_{n,0} \ \rho_{\nu+1}(x)  \over  D_n \ \rho_\nu(x)}   \sum_{r=0}^n    D_{n,r}   r! (1+\nu)_r \left[   \sum_{k=0}^{r} {B_{k,k-1} (x) \over (2k)!}   {}_3F_2 \left(-2k,\ 1+\nu+r,\ 1+r;\ 1+\nu,\ 1;\ 1 \right)\right.$$

$$\left. +  \sum_{k=0}^{r-1} {B_{k,k} (x) \over (2k+1)!}   {}_3F_2 \left(-2k-1,\ 1+\nu+r,\ 1+r;\ 1+\nu,\ 1;\ 1 \right) \right] .$$
But the existence of  a multiple orthogonal polynomial sequence with respect to the vector of weight functions $(\rho_\nu, \rho_{\nu+1})$ implies the identity

$$ \sum_{r=0}^n    D_{n,r}   r! (1+\nu)_r \left[ \sum_{k=0}^{r} {B_{k,k-1} (x) \over (2k)!}   {}_3F_2 \left(-2k,\ 1+\nu+r,\ 1+r;\ 1+\nu,\ 1;\ 1 \right)\right.$$

$$\bigg. +  \sum_{k=0}^{r-1} {B_{k,k} (x) \over (2k+1)!}   {}_3F_2 \left(-2k-1,\ 1+\nu+r,\ 1+r;\ 1+\nu,\ 1;\ 1 \right) \bigg]\equiv 0,$$
which drives to (4.43) and completes the proof. 

\end{proof}

{\bf Remark 2}.  In a similar manner orthogonal polynomials with the weight function $\rho_{\nu+1}(x)\rho_\nu(x)$ can be investigated. We leave this topic to the interested reader. 

 Finally, in this section we establish  the generating function for  polynomials $P_n$, which is  defined as usually  by the equality

$$G(x,z) = \sum_{n=0}^\infty P_n (x) {z^n\over n!} ,\quad  x >0,\  z \in \mathbb{C},\eqno(4.45)$$
where $|z| < h_x$ and $h_x >0$ is a  convergence radius of the power series. To do this,  we employ (3.1), (4.36) and (4.39), having  the following equality from (4.45) 

$$  G(x,z) =  {1 \over \rho_\nu(x) } \sum_{n=0}^\infty  {z^n\over n!}  \sum_{k=0}^{2n} {h_{2n,k} \over k! }   {d^{k}\over dx^{k}} \left[  x^{k}  \rho_\nu(x) \right]=  {1\over \rho_\nu(x) } \sum_{n=0}^\infty  {z^n\over n!}  \sum_{k=0}^{2n}  h_{2n,k}   \sum_{j=0}^k \binom{k}{j}  {(-1)^j\over (k-j)!} \  x^j \rho_{\nu-j}(x) . $$
Meanwhile, the product  $x^j\rho_{\nu-j}(x)$ is expressed in \cite{Cous}  as follows 

$$x^{ j}  \rho_{\nu-j}(x) = x^{ j/2} r_j(2\sqrt x; \nu ) \rho_{\nu}(x) +  x^{ (j-1)/2} r_{j-1} (2\sqrt x; \nu -1) \rho_{\nu+1}(x), \quad  j \in {\mathbb N}_0,$$
where $r_{-1}(z;\nu)=0$, 

$$ x^{ j/2} r_j(2\sqrt x; \nu ) = (-1)^j \sum_{i=0}^{[j/2]}   (\nu+i-j+1)_{j-2i} (j-2i+1)_i  {x^i\over i!}.\eqno(4.46)$$
Therefore this leads to the final  expression of  the generating function for the sequence $\left(P_n\right)_{n \in \mathbb{N}_0}$, namely,

$$  G(x,z) =   \sum_{n=0}^\infty  \sum_{k=0}^{2n} \sum_{j=0}^k  \binom{k}{j} { (-1)^j \ h_{2n,k} \over n!\  (k-j)! }   x^{j/2}  r_j(2\sqrt x; \nu ) z^n $$

$$+  {\rho_{\nu+1}(x) \over \rho_\nu(x) } \sum_{n=0}^\infty \sum_{k=0}^{2n} \sum_{j=0}^k  \binom{k}{j}  { (-1)^j \  h_{2n,k} \over n!\  (k-j)! }  x^{(j-1)/2} r_{j-1}(2\sqrt x; \nu-1 ) z^n,\eqno(4.47) $$
where coefficients $h_{2n,k}$ are defined by (4.40).

\section{Note on the multiple orthogonal polynomials}

In this section we will exhibit two types of multiple orthogonal polynomials for the vector of weight functions $(\rho^2_\nu,\  \rho^2_{\nu+1}, \  \rho_\nu\rho_{\nu+1}),\ \nu > -1/2 $ over $\mathbb{R}_+$ with an additional factor $x^\alpha,\ \alpha > -1$. Precisely, we consider the type $1$ polynomials $(A^\alpha_n,\ B^\alpha_{n-1},\ C^\alpha_{n-1}),\ n \in \mathbb{N}$ of degree $n,\ n-1$, respectively, satisfying the orthogonality conditions

$$\int_0^\infty \left[ A^\alpha_n(x) \rho^2_\nu (x) +  B^\alpha_{n-1}(x) \rho^2_{\nu+1} (x)+  C^\alpha_{n-1}(x) \rho_\nu (x) \rho_{\nu+1} (x) \right] x^{\alpha + m} dx = 0, \ m= 0,1,\dots, 3n-1.\eqno(5.1)$$
So, we have $3n$ linear homogeneous equations with $3n+1$ unknown coefficients of polynomials $A^\alpha_n,\ B^\alpha_{n-1},\ C^\alpha_{n-1}$. Therefore we can find type $1$ polynomials up to a multiplicative factor.  Let us denote the function $q^\alpha_{n,n-1,n-1}$ for the convenience 

$$ q^\alpha_{n,n-1,n-1}(x) =  A^\alpha_n(x) \rho^2_\nu (x) +  B^\alpha_{n-1}(x) \rho^2_{\nu+1} (x)+  C^\alpha_{n-1}(x) \rho_\nu (x) \rho_{\nu+1} (x) .\eqno(5.2)$$
Type $2$ polynomials $p^\alpha_{n,n-1,n}$ are monic polynomials of degree $3n-1$ which satisfy the multiple orthogonality conditions

$$ \int_0^\infty  p^\alpha_{n,n-1,n}(x) \rho^2_\nu (x) x^{\alpha+m} dx = 0,\quad m= 0,1,\dots, n-1,\eqno(5.3)$$

 $$ \int_0^\infty  p^\alpha_{n,n-1,n}(x) \rho^2_{\nu+1} (x) x^{\alpha+m} dx = 0,\quad m= 0,1,\dots, n-2,\eqno(5.4)$$

$$ \int_0^\infty  p^\alpha_{n,n-1,n}(x) \rho_\nu (x) \rho_{\nu+1}(x) x^{\alpha+m} dx = 0,\quad m= 0,1,\dots, n-1.\eqno(5.5)$$
This gives $3n-1$ linear equations with $3n-1$ unknown coefficients  since the leading coefficient is equal to $1$. Hence it can be uniquely determined. 

When $|\nu| < 1/2$ the uniqueness of the representation (5.2)  is validated by the following 

{\bf Theorem 4}.  {\it Let $n, m, l  \in \mathbb{N}_0, |\nu| < 1/2, \   f_n,\ g_{m},\ h_l$ be polynomials of  degree at most $n,\ m,\ l$, respectively.  Let 
 
 $$f_n(x) \rho^2_\nu(x)+ g_{m} (x) \rho^2_{\nu+1} (x) +  h_{l} (x) \rho_\nu(x) \rho_{\nu+1} (x) = 0\eqno(5.6)$$
 for all $x >0$. Then $f_n \equiv 0, \ g_{m} \equiv 0,\ h_l \equiv 0.$ }
 
 \begin{proof} As is known (cf. \cite{YAP})  the quotient $\rho_\nu/ \rho_{\nu+1}$ is represented by the Ismail integral 
 
$$ {\rho_{\nu}(x) \over \rho_{\nu+1}(x) } = {1\over \pi^2} \int_0^\infty { s^{-1} ds \over (x+s)( J_{\nu+1}^2(2\sqrt s) + 
Y_{\nu+1}^2(2\sqrt s) )},\eqno(5.7)$$
where $J_\nu, Y_\nu$ are Bessel functions of the first and second kind, respectively \cite{PrudnikovMarichev}, Vol. II.  In fact,   let $r \ge \max\{ n, m+1, l+1\}. $  Hence, dividing (5.6) by $\rho_\nu(x) \rho_{\nu+1} (x)$ and using (3.6), we find

$$f_n(x) {\rho_\nu(x)\over \rho_{\nu+1} (x)}  +  x g_{m} (x) {\rho_{\nu-1} (x)\over \rho_{\nu} (x) } + \nu g_{m} (x) + h_{l} (x) = 0.\eqno(5.8)$$
Then, differentiating $r$ times, it gives 

 $${d^{r}\over dx^{r} } \left[ f_n(x) {\rho_\nu(x)\over \rho_{\nu+1} (x)} \right] +  {d^{r}\over dx^{r} } \left[ x g_m(x) {\rho_{\nu-1} (x)\over \rho_{\nu} (x)} \right] = 0.\eqno(5.9)$$
Assuming $f_n(x)= \sum_{k=0}^n f_{n,k}\  x^k$ and employing (5.7), the first term on the left-hand side of (5.9) can be treated as follows

$$ {d^{r}\over dx^{r} } \left[ f_n(x) {\rho_\nu(x)\over \rho_{\nu+1} (x)} \right]  = {1\over \pi^2}  {d^{r}\over dx^{r} }   \sum_{k=0}^n f_{n,k} x^k \int_0^\infty e^{-xy } dy \int_0^\infty { e^{-sy} \ s^{-1} ds \over  J_{\nu+1}^2(2\sqrt s) + 
Y_{\nu+1}^2(2\sqrt s)}$$

$$=  {1\over \pi^2}    \sum_{k=0}^n f_{n,k} (-1)^k  {d^{r}\over dx^{r} }  \int_0^\infty  {d^k\over dy^k} \left[ e^{-xy } \right] dy \int_0^\infty { e^{-sy} \ s^{-1} ds \over  J_{\nu+1}^2(2\sqrt s) + Y_{\nu+1}^2(2\sqrt s)}$$

$$=  {1\over \pi^2}    \sum_{k=0}^n f_{n,k} (-1)^k   \int_0^\infty  {\partial^{k+r} \over \partial y^k \partial x^{r} } \left[ e^{-xy } \right] dy \int_0^\infty { e^{-sy} \ s^{-1} ds \over  J_{\nu+1}^2(2\sqrt s) + Y_{\nu+1}^2(2\sqrt s)}$$

 $$=  {1\over \pi^2}    \sum_{k=0}^n f_{n,k} (-1)^{k+r}    \int_0^\infty  {d^k\over dy^k} \left[  y^{r}  e^{-xy } \right] dy \int_0^\infty { e^{-sy} \ s^{-1} ds \over  J_{\nu+1}^2(2\sqrt s) + Y_{\nu+1}^2(2\sqrt s)},$$
where the differentiation  under the integral sign is possible via the absolute and uniform convergence.  Now, we integrate  $k$ times by parts in the outer integral with respect to $y$ on the right-hand side of the latter equality, eliminating the  integrated terms due to the choice of $r$,  and then differentiate  under the  integral sign in the inner integral with respect to $s$ owing to the same arguments, to obtain

$$  {d^{r}\over dx^{r} } \left[ f_n(x) {\rho_\nu(x)\over \rho_{\nu+1} (x)} \right]   =  {1\over \pi^2}     \int_0^\infty   y^{r}  e^{-xy }  \int_0^\infty { e^{-sy} \ s^{-1}  \over  J_{\nu+1}^2(2\sqrt s) + Y_{\nu+1}^2(2\sqrt s)}   \left(    \sum_{k=0}^n f_{n,k} (-1)^{k+r}  s^k \right) ds .\eqno(5.10)$$
In the same fashion the second term in (5.9) is worked out to find  $( g_m(x)=  \sum_{k=0}^m g_{m,k}\  x^k$ )

$$  {d^{r}\over dx^{r} } \left[ x g_m(x) {\rho_{\nu-1} (x)\over \rho_{\nu} (x)} \right]   =  {1\over \pi^2}     \int_0^\infty   y^{r}  e^{-xy }  \int_0^\infty { e^{-sy}  \over  J_{\nu}^2(2\sqrt s) + Y_{\nu}^2(2\sqrt s)}   \left(    \sum_{k=0}^m g_{m,k} (-1)^{k+r+1}  s^k \right) ds.\eqno(5.11)$$
Substituting (5.10), (5.11) into (5.9) and cancelling twice  the Laplace transform  via its injectivity for integrable and continuous functions \cite{Tit},  we arrive at the equality

$$x g_m(-x) \left[  J_{\nu+1}^2(2\sqrt x) + Y_{\nu+1}^2(2\sqrt x) \right] + f_n(-x) \left[   J_{\nu}^2(2\sqrt x) + Y_{\nu}^2(2\sqrt x)\right] = 0,\ x > 0.\eqno(5.12) $$
The sum of squares of Bessel functions in brackets is called the Nicholson kernel, which has the Mellin-Barnes representation by virtue of Entry 8.4.20.35 in \cite{PrudnikovMarichev}, Vol. III

$$x^k  \left[   J_{\nu}^2(2\sqrt x) + Y_{\nu}^2(2\sqrt x)\right]  = {2^{1-2k}  \cos(\pi\nu)\over 2\pi^{7/2} i} \int_{\gamma-i\infty}^{\gamma+i\infty}
\Gamma(s+k) \Gamma(s+k+\nu) \Gamma(s+k-\nu) $$

$$\times \Gamma\left({1\over 2} -s-k\right) (4x)^{-s} ds,\quad |\nu | - k < \gamma < {1\over 2} - k.\eqno(5.13)$$
Then, using on the right-hand side of (5.13) the reflection formula for gamma function  \cite{Bateman}, Vol. I,  it can be written as follows

$$ {2^{1-2k}  \cos(\pi\nu)\over 2\pi^{7/2} i} \int_{\gamma-i\infty}^{\gamma+i\infty} \Gamma(s+k) \Gamma(s+k+\nu) \Gamma(s+k-\nu) \Gamma\left({1\over 2} -s-k\right) (4x)^{-s} ds$$

$$=  {2^{1-2k}  (-1)^k \cos(\pi\nu)\over 2\pi^{7/2} i} \int_{\gamma-i\infty}^{\gamma+i\infty} \Gamma(s+k) \Gamma(s+k+\nu) \Gamma(s+k-\nu) { \Gamma(1/2 -s ) \over (1/2+s)_k } (4x)^{-s} ds.\eqno(5.14)$$
Our goal now is to shift the  contour to the right to make integration along  the straight line with $|\nu| < {\rm Re} s < 1/2$. To do this we should take into account the residues at  $k$ simple poles $s_m= -1/2- m,\  m=0, 1,\dots, k-1$ which have the values

$$\hbox{Res}_{s= s_m}  \left(\frac{\Gamma(s+k) \Gamma(s+k+\nu) \Gamma(s+k-\nu)  \Gamma(1/2 -s )} {(1/2+s) (3/2+s)\dots (s +k-1/2) } (4x)^{-s}\right)$$

$$=  \frac{\Gamma(s_m+k) \Gamma(s_m+k+\nu) \Gamma(s_m+k-\nu)  \Gamma(1/2 -s_m )} {(1/2+s_m) (1/2+s_m+1)\dots (1/2+s_m +m-1) (1/2+s_m +m+1)\dots(1/2+s_m +k-1) } (4x)^{-s_m}$$

$$= { (-1)^m (4x)^{1/2+m} \over (k-m-1)!} \ \Gamma(k-m-1/2) \Gamma(k-m-1/2+\nu) \Gamma(k-m-1/2-\nu) . $$
Therefore we get from (5.14) 

$$ {2^{1-2k}  \cos(\pi\nu)\over 2\pi^{7/2} i} \int_{\gamma-i\infty}^{\gamma+i\infty} \Gamma(s+k) \Gamma(s+k+\nu) \Gamma(s+k-\nu) \Gamma\left({1\over 2} -s-k\right) (4x)^{-s} ds$$

$$=  {2^{1-2k}  (-1)^k \cos(\pi\nu)\over 2\pi^{7/2} i} \int_{\mu-i\infty}^{\mu+i\infty} \Gamma(s+k) \Gamma(s+k+\nu) \Gamma(s+k-\nu) { \Gamma(1/2 -s ) \over (1/2+s)_k } (4x)^{-s} ds$$

$$-  {4^{1-k}  \sqrt x \cos(\pi\nu)\over \pi^{5/2}}\sum_{m=0}^{k-1} {(-1)^{k+m} (4x)^m \over (k-m-1)!}  \ \Gamma(k-m-1/2) \Gamma(k-m-1/2+\nu) \Gamma(k-m-1/2-\nu),\eqno(5.15)$$
where $|\nu| < \mu  < 1/2$.  Now, recalling  Parseval's equality for the Mellin transform and Entries 8.4.2.5, 8.4.23.27 in \cite{PrudnikovMarichev}, Vol. III,  we derive from (5.13), (5.15)

$$x^k  \left[   J_{\nu}^2(2\sqrt x) + Y_{\nu}^2(2\sqrt x)\right]  = {2^{3- 2k} (-1)^k  \sqrt x \over \pi^{3}}  \cos(\pi\nu) \int_0^\infty K^2_\nu(\sqrt t) \ {t^{k-1/2}\over 4x +t} dt$$

$$- \ {4^{1-k}  \sqrt x \cos(\pi\nu)\over \pi^{5/2}}\sum_{m=0}^{k-1} {(-1)^{k+m} (4x)^m \over (k-m-1)!}  \ \Gamma(k-m-1/2) \Gamma(k-m-1/2+\nu) \Gamma(k-m-1/2-\nu).\eqno(5.16)$$
Analogously, we find  

 $$x^{k+1}  \left[   J_{\nu+1}^2(2\sqrt x) + Y_{\nu+1}^2(2\sqrt x)\right]  = {2^{1- 2k} (-1)^k  \sqrt x \over \pi^{3}}  \cos(\pi\nu) \int_0^\infty K^2_{\nu+1}(\sqrt t) \ {t^{k+1/2}\over 4x +t} dt$$
 
 $$ - \ {4^{-k}  \sqrt x \cos(\pi\nu)\over \pi^{5/2}}\sum_{m=0}^{k} {(-1)^{k+m+1} (4x)^m \over (k-m)!}  \ \Gamma(k-m+1/2) \Gamma(k-m+1/2+\nu) \Gamma(k-m+1/2-\nu).\eqno(5.17)$$
Substituting  expressions (5.16), (5.17) in (5.12), we derive after straightforward simplifications 

$$  \int_0^\infty K^2_{\nu+1}(2 \sqrt t) \  g_m (t)  {\sqrt t \over x +t} dt +    \int_0^\infty K^2_{\nu}(2 \sqrt t) \  f_n(t)  {dt  \over \sqrt t (x +t)}$$

 $$ + \sqrt\pi \sum_{k=0}^m   \sum_{j=0}^{k} { (-x)^j \over  4^{k-j} (k-j)!}  \ \Gamma(k-j+1/2) \Gamma(k-j+1/2+\nu) \Gamma(k-j+1/2-\nu)$$

$$- \sqrt\pi \sum_{k=0}^{n-1}  \sum_{j=0}^{k} { (-x)^j \over 4^{k-j}  (k-j)!}  \ \Gamma(k-j+1/2) \Gamma(k-j+1/2+\nu) \Gamma(k-j+1/2-\nu) = 0,\  x >0.\eqno(5.18)$$
Last two terms in (5.18) are polynomials of degree $m, n-1$, respectively. Hence, differentiating through $r_1 \ge \max\{n, m+1\}$ times by $x$, we obtain

$$ {d^{r_1}\over dx^{r_1}}  \int_0^\infty K^2_{\nu+1}(2 \sqrt t) \  g_m (t)  {\sqrt t \over x +t} dt +     {d^{r_1}\over dx^{r_1}}  \int_0^\infty K^2_{\nu}(2 \sqrt t) \  f_n(t)  {dt  \over \sqrt t (x +t)} = 0.\eqno(5.19)$$
The left-hand side of (5.19) is the $r_1$-th derivative of the sum of two Stieltjes transforms which are, in turn,  two fold Laplace transforms. Consequently, fulfilling the differentiation under the integral sign in (5.19) as above owing to the absolute and uniform convergence by $x \ge x_0 >0 $,  we cancel Laplace transforms of integrable functions via the injectivity.  Then with (2.2)  it yields the equality

$$f_n(x) \rho^2_\nu(x)+ g_{m} (x) \rho^2_{\nu+1} (x) = 0,\quad x >0.\eqno(5.20)$$
Comparing with (5.6), we see that $  h_l \equiv 0.$  Further,  identity (5.20) implies that $f_n, \ g_m$ have the same positive roots, if any.    Let $x >0$ be not a root of $g_m$.  Then dividing (5.20) by $g_m$ and making a differentiation, we get

$$- 2 \rho_\nu (x) \rho_{\nu-1} (x) {f_n(x)\over g_m(x)} +   \rho^2_\nu(x) {f^\prime_n(x) g_m(x)- f_n(x) g_m^\prime(x) \over g_m^2(x)} - 2 \rho_{\nu+1}(x)\rho_\nu(x) =0.$$
Since $\rho_\nu(x) >0$, we divide the previous equation by $\rho_\nu$, multiply by $x,  g_m^2$ and employ  (3.6) to find

$$- 2 g_m(x) \rho_{\nu+1}(x) \left( f_n(x) +x g_m(x) \right) + \rho_\nu(x) \left( 2\nu f_n(x) g_m(x) +   x \left( f^\prime_n(x) g_m(x)- f_n(x) g_m^\prime(x) \right) \right) =0.\eqno(5.21)$$
However, the existence of the type 1 multiple orthogonal polynomials with respect to the vector of weight functions $(\rho_\nu, \rho_{\nu+1})$ suggests the equalities

$$ g_m(x) \left( f_n(x) +x g_m(x) \right) \equiv 0,\quad\quad  2\nu f_n(x) g_m(x) +   x \left( f^\prime_n(x) g_m(x)- f_n(x) g_m^\prime(x)\right)\equiv 0.\eqno(5.22)$$
So, if $g_m \equiv 0$, it proves the theorem. Otherwise  $f_n(x)+ x g_m(x) \equiv 0$,  and with the second equation in (5.22) we have $ x(2\nu+1) g^2_m(x) \equiv 0.$ Thus $g_m \equiv 0$ and $f_n \equiv 0$ from (5.20). Theorem 4 is proved. 

  \end{proof}
  
  {\bf Remark 3}.  The choice of $\nu$ is important. For instance, for $\nu= -1/2$  the theorem  fails.  This can be seen,  taking $f_n(x) \equiv -x,\  g_m(x)\equiv  1, \ h_l(x) \equiv 0$.

{\bf Theorem 5}. {\it Let $\nu \in [0,1/2)$. For every $\alpha > 0$

$${d\over dx} \left[ x^\alpha  q^\alpha_{n,n-1,n-1}(x) \right] = x^{\alpha-1} q^{\alpha-1}_{n,n-1,n} (x)\eqno(5.23)$$
and the following differential recurrence relations hold}

$$ A^{\alpha-1} _n(x) =   (\alpha +2\nu) A^\alpha_n(x) + x [ A^\alpha_n(x) ]^\prime   - x  C^\alpha_{n-1}(x) ,\eqno(5.24)$$ 

$$ B^{\alpha-1} _{n-1}(x) = \alpha B^\alpha_{n-1}(x)  + x  [B^\alpha_{n-1}(x) ]^\prime  - C^\alpha_{n-1}(x),\eqno(5.25)$$

$$ C^{\alpha-1} _n(x) = (\alpha+\nu)  C^\alpha_{n-1}(x) + x  [C^\alpha_{n-1}(x)]^\prime   - 2  A^\alpha_n(x)  -2x B^\alpha_{n-1}(x).\eqno(5.26)$$

\begin{proof}  From (5.1), (5.2) and integration by parts we get

$$ \int_0^\infty     q^\alpha_{n,n-1,n-1}(x) x^{\alpha+m} dx =  -  {1\over m+1} \int_0^\infty  {d\over dx} \left[  x^\alpha q^\alpha_{n,n-1,n-1}(x) \right] x^{m+1} dx = 0,$$
where the integrated terms vanish for every $\alpha > -1$ via asymptotic behavior  (2.4), (2.5).  Hence (5.1) suggests the equality 

$$\int_0^\infty  {d\over dx} \left[  x^\alpha q^\alpha_{n,n-1,n-1}(x) \right] x^{m} dx = 0, \quad    m= 1,\dots, 3n.$$
But,  evidently, 

$$\int_0^\infty  {d\over dx} \left[  x^\alpha q^\alpha_{n,n-1,n-1}(x) \right]  dx = 0,\quad \alpha > 0.$$
Therefore

$$\int_0^\infty  {d\over dx} \left[  x^\alpha q^\alpha_{n,n-1,n-1}(x) \right] x^{m} dx = 0, \quad    m= 0,\dots, 3n.\eqno(5.27)$$
Now, working out the differentiation in (5.27), involving (3.5), (3.6), we find

$$ {d\over dx} \left[  x^\alpha q^\alpha_{n,n-1,n-1}(x) \right] = \alpha\  x^{\alpha-1}  q^\alpha_{n,n-1,n-1}(x) + x^{\alpha-1}  \left(  x [ A^\alpha_n(x) ]^\prime  \rho^2_\nu (x) + x  [B^\alpha_{n-1}(x) ]^\prime \rho^2_{\nu+1} (x)+  x  [C^\alpha_{n-1}(x)]^\prime \rho_\nu (x) \rho_{\nu+1} (x) \right. $$

$$\left. + 2  A^\alpha_n(x) \rho_\nu (x) (\nu \rho_\nu(x)- \rho_{\nu+1} (x) )  -2x B^\alpha_{n-1}(x) \rho_{\nu+1} (x)\rho_\nu(x)\right. $$

$$\left. - x  C^\alpha_{n-1}(x) \rho^2_\nu (x) +   C^\alpha_{n-1}(x)  \rho_{\nu+1} (x)  (\nu \rho_\nu(x)- \rho_{\nu+1} (x) ) \right).$$
Thus

$$  {d\over dx} \left[  x^\alpha q^\alpha_{n,n-1,n-1}(x) \right] = x^{\alpha-1} \left[ A^{\alpha-1}_n(x) \rho^2_\nu (x) +  B^{\alpha-1}_{n-1}(x) \rho^2_{\nu+1} (x)+  C^{\alpha-1}_{n}(x) \rho_\nu (x) \rho_{\nu+1} (x) \right],\eqno(5.28)$$
where $A^{\alpha-1}_n(x) ,\   B^{\alpha-1}_{n-1}(x),  \ C^{\alpha-1}_{n}(x) $ are polynomials of degree at most $n, n-1, n$, respectively, being defined by formulas (5.24), (5.25), (5.26).  The linear homogeneous system (5.27) of $3n+1$ equations contains $3n+ 2$ unknowns.  Therefore up to a constant, choosing to be one,  the left-hand side of (5.28) is equal to $ x^{\alpha-1}  q^{\alpha-1}_{n,n-1,n}(x)$, and the representation (5.28) is unique by virtue of Theorem 4.  This proves (5.23) and completes the proof of Theorem 5.

\end{proof}

{\bf Remark 4}. The same analysis for the function $q^\alpha_{n,n-1,n}(x)$ does not work. In fact, in this case we derive analogously

$$\int_0^\infty  {d\over dx} \left[  x^\alpha q^\alpha_{n,n-1,n}(x) \right] x^{m} dx = 0, \quad    m= 0,\dots, 3n+1.$$
However, working out the differentiation, we will get polynomials of degree at most $n+1, n, n$, respectively, which implies $3n+4$ unknowns for $3n+2$ equations (the so-called quasi multiple orthogonal case.) 

Finally, we establish the differentiation property for the type $2$ multiple orthogonal polynomials  $p^\alpha_{n,n-1,n}$ (or $3$-orthogonal polynomials).

{\bf Theorem 6.} {\it For every $\nu \ge 0,\ \alpha > -1$ 

$${d\over dx} \left[ p^\alpha_{n,n-1,n} (x)\right] =   (3n-1) p^{\alpha+1}_{n,n-1,n-1} (x).\eqno(5.29)$$

\begin{proof} Recalling (3.5), (3.6) and asymptotic behavior of the weight functions (2.4), (2.5), we integrate by parts in (5.3), eliminating the integrated terms,  to  deduce

$$ \int_0^\infty {d\over dx} \left[  p^\alpha_{n,n-1,n}(x) \right]  \rho^2_\nu (x) x^{\alpha+1+m} dx = 0,\quad m= 0,1,\dots, n-1.\eqno(5.30)$$
Concerning equalities (5.4), (5.5),  it corresponds the following ones 

 $$ \int_0^\infty  {d\over dx} \left[  p^\alpha_{n,n-1,n}(x) \right] \rho^2_{\nu+1} (x) x^{\alpha+1+m} dx = 0,\quad m= 0,1,\dots, n-2,\eqno(5.31)$$

$$ \int_0^\infty  {d\over dx} \left[  p^\alpha_{n,n-1,n}(x) \right] \rho_\nu (x) \rho_{\nu+1}(x) x^{\alpha+1+m} dx = 0,\quad m= 0,1,\dots, n-2.\eqno(5.32)$$
 Now   $\left[  p^\alpha_{n,n-1,n}(x) \right]^\prime$ is a polynomial of degree $3n-2$ with leading coefficient $3n-1$ and by equalities (5.30), (5.31), (5.32) it satisfies orthogonality conditions (5.3), (5.4), (5.5) for the type $2$ multiple orthogonal polynomial $ p^{\alpha+1}_{n,n-1,n-1} (x)$.  Hence we get (5.29) by unicity. Theorem 6 is proved. 

\end{proof}

\bibliographystyle{amsplain}

\section{Addendum}

\title{Corrigendum to "A method of  composition orthogonality and  new sequences of orthogonal  polynomials and functions for non-classical weights"
[J. Math. Anal. Appl. 499(2021) 125032]. DOI of original article:  https://doi.org/10.1016/j.jmaa.2021.125032}}


\vspace{1cm}

In the last section of the article one should consider the type $2$ multiple monic orthogonal polynomials $p^\alpha_{n+1,n,n+1}$  of degree $3n+2$ which satisfy the  orthogonality conditions

$$ \int_0^\infty  p^\alpha_{n+1,n,n+1}(x) \rho^2_\nu (x) x^{\alpha+m} dx = 0,\quad m= 0,1,\dots, n,\eqno(1)$$

 $$ \int_0^\infty  p^\alpha_{n+1,n,n+1}(x) \rho^2_{\nu+1} (x) x^{\alpha+m} dx = 0,\quad m= 0,1,\dots, n-1,\eqno(2)$$

$$ \int_0^\infty  p^\alpha_{n+1,n,n+1}(x) \rho_\nu (x) \rho_{\nu+1}(x) x^{\alpha+m} dx = 0,\quad m= 0,1,\dots, n.\eqno(3)$$
This gives $3n+2$ linear equations with $3n+2$ unknown coefficients  since the leading coefficient is equal to $1$. Hence it can be uniquely determined.   Theorem 6 reads as follows.

{\bf Theorem 6.} {\it For every $\nu \ge 0,\ \alpha > -1$}

$${d\over dx} \left[ p^\alpha_{n+1,n,n+1} (x)\right] =   (3n+2) p^{\alpha+1}_{n+1,n,n} (x).\eqno(4)$$

\begin{proof} In fact, recalling (3.5), (3.6) and asymptotic behavior of the weight functions (2.4), (2.5), we integrate by parts in (1), eliminating the integrated terms,  to  deduce

$$ \int_0^\infty {d\over dx} \left[  p^\alpha_{n+1,n,n+1}(x) \right]  \rho^2_\nu (x) x^{\alpha+1+m} dx = 0,\quad m= 0,1,\dots, n.$$
Concerning equalities (2), (3),  it corresponds the following ones 

 $$ \int_0^\infty  {d\over dx} \left[  p^\alpha_{n+1,n,n+1}(x) \right] \rho^2_{\nu+1} (x) x^{\alpha+1+m} dx = 0,\quad m= 0,1,\dots, n-1,$$

$$ \int_0^\infty  {d\over dx} \left[  p^\alpha_{n+1,n,n+1}(x) \right] \rho_\nu (x) \rho_{\nu+1}(x) x^{\alpha+1+m} dx = 0,\quad m= 0,1,\dots, n-1.$$
 Now   $\left[  p^\alpha_{n+1,n,n+1} \right]^\prime$ is a polynomial of degree $3n+1$ with leading coefficient $3n+2$ and by the latter equalities it satisfies orthogonality conditions (1), (2), (3) for the type $2$ multiple orthogonal polynomial $ p^{\alpha+1}_{n+1,n,n} $. Hence we get (4) by unicity and complete the proof. 
 
 \end{proof}

The author would like to apologise for any inconvenience caused.

\end{document}